\documentclass[12pt]{iopart}
\usepackage{iopams}

\begin{document}

\title[Geometrisation of Chaplygin's reducing multiplier theorem]{Geometrisation of Chaplygin's reducing multiplier theorem}

\author{A V Bolsinov$^1$, A V Borisov$^2$, I S Mamaev$^{3}$}

\address{$^1$Dept. of Mathematical Sciences, Loughborough University, Loughborough,\\
 LE11 3TU,  UK}

\address{$^{2,3}$Udmurt State University, Izhevsk, Russia}

\address{$^{3}$A.\,A.\,Blagonravov Mechanical Engineering Research Institute of RAS, Moscow, Russia}
\eads{$^1$\mailto{A.Bolsinov@lboro.ac.uk}, $^2$\mailto{borisov@rcd.ru}, $^3$\mailto{mamaev@rcd.ru}}

\begin{abstract}
We develop the reducing multiplier theory for a special class of nonholonomic dynamical systems and show that the non-linear Poisson brackets naturally obtained in the framework of this approach are all isomorphic to the Lie-Poisson $e(3)$-bracket.  As two model examples, we consider the Chaplygin ball problem on the plane and the Veselova system. In particular, we obtain an integrable gyrostatic generalisation of the Veselova system.
\end{abstract}



\newcommand{\sgn}{\mathop{\rm sgn}\nolimits}
\newcommand{\arccot}{\mathop{\rm arccot}\nolimits}
\newcommand{\sparal}{\text{\scriptsize ||}}
\newtheorem{theorem}{Theorem}
\newtheorem{proposition}{Proposition}

\section*{Introduction}

In~\cite{b21-3} S.\,A.\,Chaplygin found a special class of systems with
two degrees of freedom which can be reduced to a Lagrangian and thus
Hamiltonian form by a suitable change of time $dt=\rho\,d\tau$, where $\rho$~is a
reducing multiplier depending on the coordinates. As an illustration, he
considered the problem of motion of the so-called Chaplygin sleigh, which
can be integrated by the Hamilton\,--\,Jacobi method using the
\textit{reducing multiplier method} proposed by himself.
Afterwards it was shown that a number of
systems in nonholonomic mechanics can also be represented in the form of
\textit{Chaplygin systems} or \textit{generalised Chaplygin
systems}~\cite{bm}, and thereby are conformally Hamiltonian \cite{b21-1, bm,b21-42,bfm, bmm}. Thus, the
reducing multiplier method is one of the most effective methods for
explicit \textit{Hamiltonisation} of dynamical systems.

From today's perspective, the reducing multiplier theory is a method for
finding one of the most important \textit{tensor invariants}~\cite{kozlov}
of a dynamical system~--- the \textit{Poisson structure}~\cite{bm}. At the
same time, the application of this method requires rewriting the
equations of motion in local coordinates, which usually involves extremely
cumbersome calculations. In this paper we develop the Chaplygin method for
one class of systems frequently discussed in nonholonomic mechanics, which
allows to achieve their Hamiltonisation in a much simpler way. We
shall not dwell here on the derivation of equations of motion for
nonholonomic mechanics. A fairly detailed treatment of this can be found
in~\cite{bm_kach}.

\section{Generalised Chaplygin systems}

We recall that according to~\cite{bm}, a \textit{generalised Chaplygin
system} is a mechanical system with two degrees of freedom whose equations
of motion can be written as
\begin{equation}
\label{eq21-1}
\eqalign{
&\frac d{dt}\biggl(\frac{\partial L}{\partial\dot q_1}\biggr)-\frac{\partial L}{\partial q_1}=
    \dot{q}_2S,\qquad\frac d{dt}\biggl(\frac{\partial L}{\partial\dot{q}_2}\biggr)-\frac{\partial L}{\partial q_2}=-\dot{q}_1S,\\
&S=a_1(\boldsymbol q)\dot{q}_1+a_2(\boldsymbol q)\dot{q}_2+b(\boldsymbol q),
}
\end{equation}
where $L$ is a function of generalised coordinates $\boldsymbol q=(q_1,q_2)$ and
velocities $\dot{\boldsymbol q}=(\dot q_1,\dot q_2)$, which we may call the
Lagrangian of the system. It is straightforward to verify that
this system admits an energy integral of standard form
\begin{equation}
\label{eq1_12}
E=\sum_i \frac{\partial  L}{\partial  \dot q_i}\dot q_i-L.
\end{equation}

{\bf Remark.}
A usual Chaplygin system can be obtained by a special choice of the
function $S$ (a fortiori $b(\boldsymbol q)=0$) \cite{b21-3}.  A somewhat different generalisation of the Chaplygin systems is
proposed in~\cite{b21-22,
b21-71}.
\smallskip

If there is an invariant measure with density
depending only on the coordinates, the system can be represented in
conformally Hamiltonian form \cite{bm} (for $b(\boldsymbol q)=0$ this was shown by
S.\,A.\,Chaplygin~\cite{b21-3}). To show this, we use the Legendre
transform for the initial system \eref{eq21-1}:
\[
P_i=\frac{\partial L}{\partial\dot q_i},\qquad H=\sum\limits_iP_i\dot
q_i-L\Bigr|_{\dot q_i\to P_i}.
\]
Then the equations of motion~\eref{eq21-1} can be recast as
\begin{equation}
\label{eq21-4}
\eqalign{
&\dot q_i=\frac{\partial H}{\partial P_i},\qquad\dot P_1=-\frac{\partial H}{\partial q_1}+ \frac{\partial H}{\partial P_2}S,\qquad
    \dot P_2=-\frac{\partial H}{\partial q_2}- \frac{\partial H}{\partial P_1}S,\\
&S=a_1(\boldsymbol q)\dot q_1+a_2(\boldsymbol q)\dot q_2+b(\boldsymbol q)=A_1(\boldsymbol q)P_1+A_2(\boldsymbol q)P_2+B(\boldsymbol q).
}
\end{equation}
Here $H$ coincides with the energy integral~\eref{eq1_12} expressed in terms of the new variables.

Now assume that the system admits an invariant measure with density depending only on the coordinates:
\begin{equation}
\label{eq_zv}
\mu={\mathcal N}(\boldsymbol q)\,dP_1\,dP_2\,dq_1\,dq_2.
\end{equation}
In this case the Liouville equation for ${\mathcal N}(\boldsymbol q)$ reduces to
\[
\dot q_1\left(\frac1{\mathcal N}\frac{\partial{\mathcal N}}{\partial q_1}-A_2(q)\right)+
    \dot q_2\left(\frac1{\mathcal N}\frac{\partial{\mathcal N}}{\partial q_2}+A_1(q)\right)=0,
\]
and since ${\mathcal N}$ depends only on the coordinates, each of the brackets must vanish separately:
\begin{equation}
\label{eq4_12}
\frac1{\mathcal N}\frac{\partial{\mathcal N}}{\partial q_1}-A_2(\boldsymbol q)=0,\qquad
\frac1{\mathcal N}\frac{\partial{\mathcal N}}{\partial q_2}+A_1(\boldsymbol q)=0.
\end{equation}

Let us now make the change of variables
\[
P_i=\frac{p_i}{{\mathcal N}(\boldsymbol q)},\qquad i=1,2.
\]

Denote the Hamiltonian in the new variables as $\overline{H}(\boldsymbol q,\boldsymbol p)=H\big(\boldsymbol q,\boldsymbol P(\boldsymbol q,\boldsymbol p)\big)$. Then the following relations hold for the derivatives
\[
\frac{\partial  H}{\partial  P_i}={\mathcal N}\,\frac{\partial  \overline H}{\partial  p_i},\qquad
\frac{\partial  H}{\partial  q_i}=\frac{\partial  \overline H}{\partial  q_i}+
\frac1{{\mathcal N}}\,\frac{\partial  N}{\partial  q_i}\left(\frac{\partial  \overline H}{\partial  p_1}p_1+
    \frac{\partial  \overline H}{\partial  p_2}p_2\right).
\]

Substituting them into~\eref{eq21-4} and using~\eref{eq4_12},
we obtain
\[
\fl\eqalign{
\dot q_i={\mathcal N}(\boldsymbol q)\frac{\partial \overline H}{\partial  p_i},\\
\dot p_1 ={\mathcal N}(\boldsymbol q)\left(-\frac{\partial \overline H}{\partial  q_1}+{\mathcal N}(\boldsymbol q)B(\boldsymbol q) \,\frac{\partial \overline H}{\partial  p_2}\right),\qquad
\dot p_2 ={\mathcal N}(\boldsymbol q)\left(-\frac{\partial \overline H}{\partial  q_2}-{\mathcal N}(\boldsymbol q)B(\boldsymbol q) \,\frac{\partial \overline H}{\partial  p_1}\right).
}
\]

Thus, the following result holds.

\begin{theorem}
If the system~\eref{eq21-4} admits an invariant measure of the form~\eref{eq_zv}, it can be represented in conformally Hamiltonian form
\[
\dot q_i={\mathcal N}(\boldsymbol q)\{q_i,\overline H\},\qquad \dot p_i={\mathcal N}(\boldsymbol q)\{p_i,\overline H\},\qquad i=1,2,
\]
where the Poisson brackets are given by
\[
\{q_i,p_j\}=\delta_{ij},\qquad\{q_i,q_j\}=0,\qquad\{p_1,p_2\}={\mathcal N}(\boldsymbol q)B(\boldsymbol q).
\]
\end{theorem}

{\bf Proof.} The proof is a straightforward verification of the Jacobi identity. $\blacksquare$

\section{The Chaplygin system on $\boldsymbol{T^*S^2}$}\label{sect2}

We now consider a system which is described by means of two three-dimensional vectors
$\boldsymbol M$ and $\bgamma$ and whose equations of motion are
\begin{equation}
\label{Voz13}
\dot{\boldsymbol M}=(\boldsymbol M-S\bgamma)\times \frac{\partial {H}}{\partial{\boldsymbol M}}+\bgamma\times\frac{\partial{H}}{\partial{\bgamma}},\qquad
\dot{\bgamma}=\bgamma\times\frac{\partial{H}}{\partial{\boldsymbol M}},
\end{equation}
where the ``Hamiltonian''~$H(\boldsymbol M,\bgamma)$~is an arbitrary function (quadratic and non-degenerate in
$\boldsymbol M$) and $S(\boldsymbol M,\bgamma)$ is a function linear in $\boldsymbol M$:
\[
S=\big(\boldsymbol K(\bgamma),\boldsymbol M\big)=K_1(\bgamma)M_1+K_2(\bgamma)M_2+K_3(\bgamma)M_3.
\]

It can be proved by a straightforward verification that the system~\eref{Voz13} always admits
three integrals of motion:
\[
F_1=\bgamma^2,\qquad F_2=(\boldsymbol M,\bgamma),\qquad F_3=H(\boldsymbol M,\bgamma).
\]

Without loss of generality we can set $\bgamma^2=1$, so that
equations~\eref{eq4_12} govern the dynamical system on the family of four-dimensional manifolds
\[
{\mathcal M}_c^4=\{\boldsymbol M,\bgamma\mid \bgamma^2=1,(\boldsymbol M,\bgamma)=c\},
\]
each of which is diffeomorphic to $TS^2$.

If the entire set of variables is denoted as $\boldsymbol x=(\boldsymbol M,\bgamma)$, then
equations~\eref{Voz13} can be represented in the skew-symmetric form
\begin{eqnarray}
\dot {\boldsymbol x}=\mathbf{P}_0\,\frac{\partial H}{\partial\boldsymbol x},\nonumber\\
\mathbf{P}_0=
    \left(\begin{array}{cc}
        \mathbf{M} & \bGamma\\
        \bGamma & 0
    \end{array}\right)-
    S(\boldsymbol x)
    \left(\begin{array}{cc}
        \bGamma & 0\\
        0 & 0
    \end{array}\right)\!,\nonumber\\
\mathbf{M}=
    \left(\begin{array}{ccc}
        0 & -M_3 & M_2\\
        M_3 & 0 & -M_1\\
        -M_2 & M_1 & 0
    \end{array}\right)\!,\!\qquad
    \bGamma=
    \left(\begin{array}{ccc}
        0 & -\gamma _3 & \gamma _2\\
        \gamma _3 & 0 & -\gamma _1\\
        -\gamma _2 & \gamma _1 & 0
    \end{array}\right)\!.\nonumber
\end{eqnarray}
Here the first term is a standard Poisson structure corresponding to
the Lie algebra $e(3)$. Moreover, $\mathbf{P}_0$ additionally satisfies the equations
\[
\mathbf{P}_0\,\frac{\partial  F_1}{\partial \boldsymbol x}=0,\qquad
\mathbf{P}_0\,\frac{\partial  F_2}{\partial \boldsymbol x}=0.
\]

As above,  assume that \eref{Voz13} admits an invariant measure with density depending only on $\bgamma$:
\begin{equation}
\label{eq7_12}
\mu=\rho(\bgamma)\,d\boldsymbol M\,d\bgamma.
\end{equation}
In this case the Liouville equation for the vector field $\boldsymbol V(\boldsymbol M,\bgamma)$ defined
by the system~\eref{Voz13} can be represented as
\[
\mathrm{div}\, \rho \boldsymbol V=
    \left(\frac{\partial  H}{\partial  \boldsymbol M}, \rho\bgamma\times \boldsymbol K-\bgamma\times \frac{\partial \rho}{\partial \bgamma}\right)=0.
\]

Hence, owing to non-degeneracy of the Hamiltonian in $\boldsymbol M$, we obtain
the vector equation
\begin{equation}
\label{eq_zv2}
\bigg(\frac1{\rho}\frac{\partial  \rho}{\partial \bgamma}-\boldsymbol K\bigg)\times \bgamma=0.
\end{equation}
Using this relation, we can prove by direct computation

\begin{proposition}
If  $\rho(\bgamma)$ satisfies equation~\eref{eq_zv2}, then the
tensor $\mathbf{P}=\frac1{\rho(\bgamma)}\mathbf{P}_0$ satisfies the Jacobi
identity and therefore is a Poisson structure on $\mathbb R^6(\boldsymbol M,\bgamma)$.
\end{proposition}

Thus, we finally obtain

\begin{theorem}
If the system~\eref{Voz13} admits an invariant measure \eref{eq7_12} with
density depending only on $\bgamma$, it can be represented in the conformally
Hamiltonian form
\[
\dot {\boldsymbol x}=\rho(\bgamma)\mathbf{P}(\boldsymbol x)\frac{\partial  H}{\partial  \boldsymbol x},
\]
where $\mathbf{P}(\boldsymbol x)=\rho^{-1}\mathbf{P}_0(\boldsymbol x)$~is a Poisson structure of rank 4 with
the Casimir functions
\[
F_1=\bgamma^2,\qquad
F_2=(\boldsymbol M,\bgamma).
\]
\end{theorem}

Equation~\eref{eq_zv2} can be solved for the vector $\boldsymbol K$
as follows:
\[
\boldsymbol K=\rho f(\bgamma)\bgamma + \frac{1}{\rho}\,\frac{\partial  \rho}{\partial  \bgamma},
\]
where $f(\bgamma)$~is an arbitrary function. Thus, we have naturally obtained
a special class of Poisson structures on the space
$\mathbb{R}^6(\boldsymbol M,\bgamma)$, which can be written as
\begin{equation}
\label{eq9_12}
\mathbf{P}=\frac{1}{\rho}
\left(\begin{array}{cc}
\mathbf{M} & \bGamma\\
\bGamma & 0
\end{array}\right)
-\left(\frac{1}{\rho^2}\,\frac{\partial \rho}{\partial \bgamma}+f(\bgamma)\bgamma,\boldsymbol M\right)
\left(\begin{array}{cc}
\bGamma & 0\\
0 & 0
\end{array}\right)\!.
\end{equation}

{\bf Remark.} If we add a term of the form
\[
\Phi(\bgamma)
\left(\begin{array}{cc}
\bGamma & 0\\
0 & 0
\end{array}\right)\!
\]
to the bracket~\eref{eq9_12},
where $\Phi(\bgamma)$~is an arbitrary function, then the Jacobi identity will still hold. We use such a  modification of
\eref{eq9_12} below, see \eref{newS}.

\smallskip

We give two examples.

\medskip
\textbf{The problem of the Chaplygin ball on a plane~\cite{chaplygin1}}
describing the rolling of a balanced dynamically asymmetric ball without
slipping on a horizontal plane. \smallskip

In appropriate variables the equations of motion can be represented in the
form~\eref{Voz13}, see~\cite{bts,bmg,bm_kach} with
\begin{equation}
\label{Voz12}
\fl H_1=\frac{1}{2}\bigg((\mathbf{A}\boldsymbol M, \boldsymbol M)+\frac{(\mathbf{A}\boldsymbol M,\bgamma)^2}{{\mathcal D}^{-1}-(\bgamma,\mathbf{A}\bgamma)}\bigg)+
    U_1(\bgamma), \qquad
S=\frac{(\mathbf{A}\boldsymbol M, \bgamma)}{{\mathcal D}^{-1}-(\mathbf{A}\bgamma, \bgamma)},
\end{equation}
where ${\mathcal D}={\rm const}$, $\mathbf{A}$~is a constant diagonal
matrix. The ball's angular momentum $\boldsymbol M$ relative to the point of contact
is expressed in terms of the physical variable $\bomega$, angular
velocity, by the formula
\begin{equation}
\label{eq10_12_}
\boldsymbol M=\mathbf{A}^{-1}\bomega-{\mathcal D}(\bomega,\bgamma)\bgamma,\qquad
\bomega=\mathbf{A}(\boldsymbol M + S\bgamma),
\end{equation}
and the integral $(\boldsymbol M,\bgamma)$ can take arbitrary values.

The density of the invariant measure~\eref{eq7_12} of the system and the
function $f(\bgamma)$ for the bracket~\eref{eq9_12} has the form
\[
\rho=\frac1{\displaystyle\sqrt{{\mathcal D}^{-1}-(\bgamma,\mathbf{A}\bgamma)}},\qquad   f(\bgamma)=0.
\]

\medskip
\textbf{The Veselova system~\cite{b21-5,b21-7,fedorov}} governing the
dynamics of a body with a fixed point subject to the nonholonomic
constraint $(\bomega,\bgamma)=b={\rm const}$, where $\bomega$~is the
angular velocity of the ball and $\bgamma$~is a unit vector fixed in space.
\smallskip

In the body-fixed frame, the equations of motion can be represented in the form~\eref{Voz13}, see~\cite{mamaev2} with
\begin{equation}
\label{eq10_12}
\fl H_2=\frac12\left((\boldsymbol M,\widehat{\mathbf{A}}\boldsymbol M)-
    \frac{\big((\widehat{\mathbf{A}}-\mathbf{E})\boldsymbol M,\bgamma\big)^2}{(\widehat{\mathbf{A}}\bgamma,\bgamma)}\right)+U_2(\bgamma),\qquad
S=-\frac{\big((\widehat{\mathbf{A}}-\mathbf{E})\boldsymbol M,\bgamma\big)}{(\widehat{\mathbf{A}}\bgamma,\bgamma)},
\end{equation}
where $\mathbf{A}=\mathbf{I}^{-1}$~is the constant matrix inverse to the tensor of inertia, and
the angular momentum $\boldsymbol M$ is expressed in terms of the angular velocity $\bomega$ of the body as follows
\begin{equation}
\label{eq12_}
\boldsymbol M=\widehat{\mathbf{A}}^{-1}\bomega+\big((\widehat{\mathbf{A}}^{-1}-\mathbf{E})\bomega,\bgamma)\bgamma,\qquad
 \bomega=\widehat{\mathbf{A}}(\boldsymbol M - S\bgamma),
\end{equation}
where the area integral coincides with the constraint equation:
\[
(\boldsymbol M,\bgamma)=(\bomega,\bgamma)=b.
\]

The density of the invariant measure~\eref{eq7_12} and the function $f(\bgamma)$ coincide in this case:
\[
\rho(\bgamma)=f(\bgamma)=\frac1{\displaystyle\sqrt{(\bgamma,\widehat{\mathbf{A}}\bgamma)}}.
\]

A Lagrangian representation for $b=0$ after a change of time was obtained
in~\cite{b21-1}, the corresponding conformally Hamiltonian representation
in~\cite{mamaev2}, and another conformally Hamiltonian representation was
found in~\cite{bm}.

If the potential $U(\bgamma)$ for these systems is not zero, then, as a
rule, the corresponding equations of motion turn out to be nonintegrable so that
this Hamiltonisation method is essentially different from that used
in~\cite{biz_tsig}, where the existence of a complete set of first integrals was
required.

We also note that if one makes a change of the parameters and the potential in the Chaplygin ball problem:
\[
\mathbf{A}={\mathcal D}^{-1}(\mathbf{E}-\widehat{\mathbf{A}}),\qquad
U_1(\bgamma)={\mathcal D}^{-1}U_2(\bgamma),
\]
then we find that the Hamiltonian~\eref{Voz12} becomes
\begin{equation}
\label{eq11_12}
H_1=\frac{{\mathcal D}^{-1}}{2}(\boldsymbol M,\boldsymbol M)-{\mathcal D}^{-1}H_2,
\end{equation}
and the Poisson structure of the Chaplygin ball is transformed into a
Poisson structure of the Veselova system. Consequently, these two systems
are defined on the same Poisson manifold~\cite{tsiganov}, and their
Hamiltonians are related by~\eref{eq11_12}. If $U_i(\bgamma)=0$, $i=1,2$, then the function
$F=\boldsymbol M^2$ is an integral for the both systems, which implies that their trajectories
turn out to be rectilinear windings (transverse to each other) on the same invariant 
tori~\cite{fedorov}.

\section{Reduction to the  $e(3)$-bracket}

Introducing new notation $g = \rho^{-1}$,  we can rewrite the Poisson
structure~\eref{eq9_12} in a shorter form that is more convenient for further analysis
\begin{equation}
\label{Pgf}
{\mathbf P} =  g
\left(\begin{array}{cc}
\mathbf{M} & \bGamma\\
\bGamma & 0
\end{array}\right)
+ \bigg(\frac{\partial  g}{\partial \bgamma} - f\cdot \bgamma, \boldsymbol M\bigg)
\left(\begin{array}{cc}
\bGamma & 0\\
0 & 0
\end{array}\right)\!.
\end{equation}

Let us examine the family of such Poisson structures in more detail. First of
all, we see that this family is parametrised by two arbitrary functions $g(\bgamma) >0$ and $f(\bgamma)$
and we will denote the
corresponding Poisson structures by ${\mathbf P}_{g,f}$. Notice that all ${\mathbf P}_{g,f}$
possess the same Casimir functions $(\boldsymbol M,\bgamma)$ and
$(\bgamma, \bgamma)$.

For simplicity we confine our attention to the physical case $\bgamma^2
=(\bgamma,\bgamma)=1$, that is, we restrict all the objects to the five-dimensional (Poisson)
manifold $S^2(\bgamma)\times \mathbb{R}^3(\boldsymbol M)$.

One of our goals is to find out to what canonical form these Poisson
structures can be reduced. First of all, we note that the symplectic
leaves of ${\mathbf P}_{g,f}$  are all diffeomorphic to the cotangent bundle to
the sphere $T^*S^2$.  From the explicit form~\eref{Pgf} of the Poisson
structure it may be inferred that the symplectic structure on each leaf
$T^*S^2$ will be the sum of the canonical form $dp\wedge dq$ and a magnetic
term, that is, a closed 2-form $\omega_{\mathrm{magn}}$ on the sphere. By
the Moser theorem \cite{moser}, such forms $\omega_{\mathrm{magn}}$
are parametrised up to a symplectomorphism by one single number, namely  $\displaystyle\int_{S^2}\omega_{\mathrm{magn}}$. Thus,
for each Poisson structure we have a one-parameter family of symplectic
leaves whose type is also defined by exactly one parameter.  This
observation leads us to the conjecture that by ``redistributing'', if
necessary, the symplectic leaves and then by applying a certain
symplectomorphism to each single symplectic leaf, we can transform any Poisson structure 
${\mathbf P}_{g,f}$ to any other ${\mathbf P}_{\widetilde  g,
\widetilde  f}$.

{\bf Remark.}
On the zero level $(\boldsymbol M,\bgamma)=0$, the Poisson structure~\eref{Pgf}
is reduced to the canonical $e(3)$-bracket by a very simple transformation~\cite{bm}:
\[
(\boldsymbol M,\bgamma)\mapsto \big(g^{-1}(\bgamma)\boldsymbol M,\bgamma\big).
\]
Thus, for the Chaplygin ball we have:
\[
(\boldsymbol M,\bgamma)\mapsto \Big(\big({\mathcal D}^{-1}-(\bgamma,\mathbf{A}\bgamma)\big)^{-1/2}\boldsymbol M,\bgamma\Big),
\]
and for the Veselova system:
\[
(\boldsymbol M,\bgamma)\mapsto \big((\bgamma,\mathbf{A}\bgamma)^{-1/2}\boldsymbol M,\bgamma\big).
\]
\smallskip

We start by describing a class of natural transformations which preserve
the form of ${\mathbf P}_{g,f}$,  but change the
parameters $g$ and~$f$. Consider the transformations of the form
\begin{equation}
\label{transf}
(\boldsymbol M,\bgamma)  \mapsto (\widetilde {\boldsymbol M}, \bgamma), \qquad  \widetilde {\boldsymbol M} = \mathbf{A}(\bgamma) \boldsymbol M,
\end{equation}
where $\mathbf{A}(\bgamma)$ is a linear operator in $\mathbb{R}^3$ whose components depend on $\bgamma$.

\begin{proposition}
\label{pro_2}
For each point $\bgamma\in S^2$, consider the orthogonal
decomposition $\boldsymbol M=\boldsymbol M' + \boldsymbol M''$,  where $\boldsymbol M''= (\boldsymbol M,\bgamma) \bgamma$ is
the projection of $\boldsymbol M$ onto the vector $\bgamma$, and $\boldsymbol M' = \boldsymbol M- \boldsymbol M''$
is the projection of $\boldsymbol M$ onto the plane perpendicular to $\bgamma$, i.e., the tangent plane $T_\gamma S^2$. Let
$$
\widetilde  {\boldsymbol M }=  \alpha (\bgamma) \boldsymbol M' + c \boldsymbol M'' + \boldsymbol M'' \times \boldsymbol h(\bgamma),
$$
where $c\ne 0$ is a constant, $\alpha(\bgamma)>0$ is an arbitrary scalar
function, and $\boldsymbol h(\bgamma)$ is an arbitrary vector function of $\bgamma$.
Then the transformation \eref{transf} sends $\mathbf P_{g,f}$
to a Poisson structure ${\mathbf P_{\widetilde  g,
\widetilde  f}}$ of the same kind with parameters
\begin{equation}
\label{transform}
\fl \widetilde  g  = \alpha  g ,\qquad
\widetilde  f = \frac{\alpha^2}{c}  f  +
    \left(\frac{\alpha}{c} -1\right)\bigg( \widetilde  g -\bigg (\bgamma, \frac{\partial  \widetilde g}{\partial  \bgamma}\bigg)\bigg)+
    \frac{1}{c}\bigg(\bgamma,   \widetilde  g \, \frac{\partial \alpha}{\partial \bgamma}+
    \widetilde  g^2 \,  \mathrm{curl}\, \bigg(\frac{\boldsymbol h}{ \widetilde  g}\bigg) \bigg).
\end{equation}
\end{proposition}

The proof of Proposition \ref{pro_2} is a straightforward
verification and we confine ourselves to commenting on the geometric meaning of the
transformation $\boldsymbol M \mapsto \widetilde  {\boldsymbol M}=\mathbf{A}(\bgamma){\boldsymbol M}$ used in this proposition.
Consider the orthonormal basis $\boldsymbol e_1, \boldsymbol e_2, \boldsymbol e_3$ related to the vector
$\bgamma$ in the space $\mathbb{R}^3(\boldsymbol M)$.  Namely, $\boldsymbol e_1$ and $\boldsymbol e_2$
are two orthonormal vectors lying in a tangent plane to the unit sphere at
point $\bgamma$, and $\boldsymbol e_3$ is the normal vector to this sphere at the
same point, i.e., $e_3=\bgamma$. In this basis the matrix of 
$\mathbf{A}=\mathbf{A}(\bgamma)$ has the form
\[
\mathbf{A}=
\left(\begin{array}{ccc}
    \alpha & 0 &  a \\
    0 & \alpha & b \\
    0 & 0 & c
\end{array}\right)
\]
where $\alpha$, $a$ and $b$ depend on $\bgamma$, and $c$ is constant.

This is exactly the general form of the transformation $\mathbf{A}$ which
satisfies our requirements. Indeed, the Casimir function $(\boldsymbol M,\bgamma)$
should be mapped to itself with possible multiplication by some constant $c$.
Therefore, the plane defined by the equation $(\boldsymbol M,\bgamma)=0$ is sent to
itself,  and in the orthogonal direction the transformation is a
dilatation with ratio $c$ independent of $\bgamma$.  These
conditions completely define the last row of the matrix $\mathbf{A}$.

Furthermore, the relations $\{ M_i, \gamma_j\} = - g
\cdot\varepsilon_{ijk} \gamma_k$ can be formally rewritten in vector form
as $\{ \boldsymbol M, \bgamma\} = - g \, \boldsymbol M\times\bgamma$.  Since their form must
remain the same, we obtain the condition
\[
g \, \bigl(\mathbf{A}(\bgamma) \boldsymbol M\bigr)\times  \bgamma =  \widetilde  g \, \boldsymbol M\times\bgamma.
\]
This means that on the tangent plane $T_{\bgamma}S^2$ the operator $\mathbf{A}$
must act as multiplication by some number $\alpha$ (depending on
$\bgamma$). There are no restrictions on the elements $a$ and $b$,  they
are given by the vector function~$\boldsymbol h$  (this function itself has
3~components, but only two of them are significant, since nothing is changed
by adding to $\boldsymbol h$ any vector proportional to $\bgamma$).

\medskip

Notice that the set of
transformations described in Proposition~\ref{pro_2} forms a group (which
is, of course, infinite-dimensional, since its parameters contain
arbitrary functions $\alpha$ and $\boldsymbol h$).  It is easily verified that
performing successively two transformations with parameters $(\alpha_1, c_1, \boldsymbol h_1)$ and
$(\alpha_2, c_2, \boldsymbol h_2)$ is equivalent to the transformation with
parameters $(\alpha_1\alpha_2, c_1c_2,  \boldsymbol h_1\alpha_2 + \boldsymbol h_2 c_1)$. The
above-mentioned rule specifies a group binary operation, which simply copies the matrix
multiplication:
\[
\left(\begin{array}{cc}
    \alpha_2 & \boldsymbol h_2 \\
    0 & c_2
\end{array}\right)
\left(\begin{array}{cc}
    \alpha_1 & \boldsymbol h_1 \\
    0 & c_1
\end{array}\right) =
\left(\begin{array}{cc}
    \alpha_1 \alpha_2& \boldsymbol h_1\alpha_2  + \boldsymbol h_2 c_1  \\
    0 & c_1 c_2
\end{array}\right)
\]
This group acts in a natural way on the family of
Poisson structures $\left\{ {\mathbf P}_{g,f} \right\}$ or, which is the same, on the space of parameters $g,
f$. The above relations \eref{transform} can be understood as explicit formulae for this action.  If the action is formally denoted by $(\widetilde  g,
\widetilde  f) = \Psi_{(\alpha, c, \boldsymbol h)} (g, f)$, then, as is easily verified by successively performing two transformations, it satisfies the standard action rule. Namely, if
\[
(\widetilde  g, \widetilde  f) = \Psi_{(\alpha_1, c_1, \boldsymbol h_1)} (g, f) \qquad \mbox{and} \qquad
(\widetilde {\widetilde  g}, \widetilde {\widetilde  f}) = \Psi_{(\alpha_2, c_2, \boldsymbol h_2)} (\widetilde  g, \widetilde  f),
\]
then
\[
(\widetilde {\widetilde  g}, \widetilde {\widetilde  f}) = \Psi_{(\alpha_1\alpha_2, c_1c_2,  \boldsymbol h_1\alpha_2 + \boldsymbol h_2 c_1)} (g,f).
\]

For an explicit verification of this fact it is convenient to rewrite \eref{transform} as
\[
\fl \widetilde  g  = \alpha g,  \qquad
\widetilde  f = \frac{\alpha^2}{c} \bigg( f  +  g - \bigg(\bgamma, \frac{\partial  g}{\partial  \bgamma} \bigg)\bigg) -
\bigg( \widetilde  g - \bigg(\bgamma, \frac{\partial \widetilde g}{\partial \bgamma}\bigg)\bigg) +
\frac{\widetilde  g^2}{c} \bigg(\bgamma, \mathrm{curl} \bigg(\frac{\boldsymbol h}{\widetilde  g}\bigg) \bigg).
\]
Now the verification presents no difficulty.

From the viewpoint of group theory it would now be natural to ask the
question: {\it what are the orbits of this action\/}?
In other words, we want to understand which Poisson structures may be
transferred to each other by the above-mentioned transformations. The answer
turns out to be very simple:  \textit{the action described above has one
single orbit}, i.e., all Poisson structures in this family
are equivalent to each other.  In particular, the following
theorem holds:

\begin{theorem}\label{thm3}
Every Poisson structure $\mathbf P_{g,f}$ of the form \eref{Pgf} on the level
$\bgamma^2 =1$ is isomorphic to the standard Lie-Poisson structure $\mathbf
P_{1,0}$ related to the Lie algebra $e(3)$.
\end{theorem}

{\bf Proof.}  It is sufficient to choose parameters $(\alpha, c, \boldsymbol h)$ in \eref{transform} in such a way that $\widetilde
g=1$ and $\widetilde  f=0$. The first condition immediately defines the
function $\alpha$, namely, $\alpha = g^{-1}$. After that the second
condition reduces to
\[
\frac{\alpha^2}{c} f + \left( \frac{\alpha}{c} - 1\right) +
\frac{1}{c}\bigg(\bgamma, \frac{\partial \alpha }{\partial \bgamma}\bigg) + \frac{1}{c} (\bgamma, \mathrm{curl}\, \boldsymbol h) = 0
\]
or, equivalently,
\[
\alpha^2 f + \alpha + \bigg(\bgamma,\frac{\partial \alpha }{\partial \bgamma}\bigg) - c +  (\bgamma, \mathrm{curl}\, \boldsymbol h) = 0,
\]
where the constant $c$ and the vector function $\boldsymbol h$ are the unknowns. This equation can now be rewritten as
\begin{equation}
\label{eqzv}
(\bgamma, \mathrm{curl}\, \boldsymbol h) =  F(\bgamma) + c,
\end{equation}
where $F(\bgamma)$ is a given function. Notice that \eref{eqzv} has to be fulfilled only on the unit sphere $\bgamma^2 = 1$. The conditions
for solving the equations of this form are well known. In the
differential-geometric sense this equation simply means that we are
looking for an antiderivative of the $2$-form $(F + c)\,d\sigma$ on the unit
sphere, where $d\sigma$~is the standard area form. Such a 1-form can be found
if and only if $\displaystyle\int_{S^2} (F+c) \,d\sigma = 0$. This
condition can always be achieved by choosing a constant $c$.  $\blacksquare$

{\bf Remark.}
In a similar manner, the bracket $\mathbf P_{g,f}$ can be reduced to the
standard form on the whole space  $\mathbb R^6({\boldsymbol M},\bgamma )\simeq e^*(3)$, i.e., without the
additional restriction $\bgamma^2 =1$. To that end, we have to extend the
class of transformations by assuming that $c$ depends on
$\bgamma^2$. Since $\bgamma^2$ is a Casimir function,  $c(\bgamma^2)$ may
be treated, as before, as a constant and hence the
formulae do not essentially change. The conditions for solvability of the
equation $(\bgamma, \mathrm{curl}\, \boldsymbol h) =  F(\bgamma) + c(\bgamma^2)$
remain the same, but now they have to be verified on the spheres of all
radii.  As before, we are able to ensure that they are satisfied, since
the necessary constants can now be chosen depending on the square $\bgamma^2$ of the radius
$|\bgamma|=\sqrt{\gamma_1^2+\gamma_2^2+\gamma_3^2}$.
\smallskip

For the Chaplygin ball problem, the function $F(\bgamma)$ in Eq.~\eref{eqzv} has the form
\[
F(\bgamma)=-\frac{{\mathcal D}^{-1}}{\big({\mathcal D}^{-1}-(\bgamma,\mathbf{A}\bgamma)\big)^{3/2}}.
\]

The solutions of Eq.~\eref{eqzv} for the unknowns $c$ and $\boldsymbol h$ can be
expressed in this case in terms of complete and incomplete elliptic
integrals. Thus, although theoretically it is not
difficult to prove reducibility of $\mathbf
P_{g,f}$ to the $e(3)$-bracket, in practice the resulting transformation can turn out
to be extremely unwieldy and non-algebraic.

\section{Generalisation to the case of a gyrostat}

In this section we consider the dynamical systems obtained by adding a rotor with constant gyroscopic
momentum $\boldsymbol k$  to the Chaplygin ball and a rigid body in
the Veselova problem. A detailed derivation of the equations of motion for
compound bodies can be found in books~\cite{dtt,wittenburg,lcha}.

The new equations with gyrostatic terms take the following form
\begin{equation}
\label{withgyr}
\dot{\boldsymbol M}=(\boldsymbol M+\boldsymbol k-S\bgamma)\times\frac{\partial H}{\partial \boldsymbol M}+
    \bgamma\times\frac{\partial H}{\partial \bgamma},\quad
\dot{\bgamma}=\bgamma\times\frac{\partial H}{\partial \boldsymbol M},
\end{equation}
where the new ``Hamiltonian'' $H$ and 
function $S$ may now depend on the gyrostatic momentum $\boldsymbol k$ as a parameter, but preserve their original structure as in Section \ref{sect2}. In particular,
\begin{equation}
\label{newS}
S=\frac{1}{g}\bigg(-\frac{\partial  g}{\partial \bgamma} + f(\bgamma)\bgamma,\boldsymbol M\bigg) + \frac{1}{g} \,\Phi(\bgamma)
\end{equation}
for some smooth functions $g(\bgamma), f(\bgamma)$ and $\Phi(\bgamma)$.

A direct calculation shows that this system remains conformally Hamiltonian. Namely, (\ref{withgyr}) can be rewritten as
$$
\dot {\boldsymbol x}=g^{-1} \mathbf{P}_{\boldsymbol k}(\boldsymbol x)\,\frac{\partial  H}{\partial  \boldsymbol x},
$$
with the Poisson structure
$\mathbf P_k$ of a more general form
\begin{equation}
\label{eq13_}
\eqalign{
\mathbf{P}_{\boldsymbol k}(\boldsymbol x)=g
    \left(\begin{array}{cc}
    \mathbf{M}_{\boldsymbol k} & \bGamma\\
    \bGamma & 0
    \end{array}\right)
    -gS
    \left(\begin{array}{cc}
    \bGamma & 0\\
    0 & 0
    \end{array}\right)\!,\\
\mathbf{M}_{\boldsymbol k}=
    \left(\begin{array}{ccc}
    0 & -M_3-k_3 & M_2+k_2\\
    M_3+k_3 & 0 & -M_1-k_1\\
    -M_2-k_2 & M_1+k_1 & 0
    \end{array}\right)\!,
}
\end{equation}
where $\boldsymbol x=(\boldsymbol M,\bgamma)$~is a complete set of variables.

The Jacobi identity for  $\mathbf{P}_{\boldsymbol k}$ is fulfilled, and the Casimir functions are
\[
F_1=\bgamma^2,\qquad   F_2=(\boldsymbol M + \boldsymbol k,\bgamma).
\]

The new expressions for $H$ and $S$ presented below can be obtained by using the methods developed in \cite{dtt,wittenburg,lcha}. We omit this computation.

For the {\it Chaplygin ball}, the vector $\boldsymbol M$ is still expressed in terms of the
angular velocity $\bomega$ by means of~\eref{eq10_12_}, and
the Hamiltonian \eref{Voz12} also remains the same. For the
bracket~\eref{eq13_} we set
\[
g=\sqrt{\displaystyle {\mathcal D}^{-1}-(\bgamma,\mathbf{A}\bgamma)},\qquad  f(\bgamma)=0, \qquad \Phi(\bgamma)=0.
\]

In other words, all the ingredients remain unchanged except for the additional terms involving ${\boldsymbol k}$ in the bracket (\ref{eq13_}). Thus, to obtain the gyrostatic generalisation of the Chaplygin ball we simply need
to replace $\mathbf{M}$ by $\mathbf{M}_{\boldsymbol k}$ in \eref{Pgf}.

For the {\it Veselova system}, when a gyrostat is added, the situation becomes less trivial and the
relations~\eref{eq12_} as well as $H_2$ and $S$ given by \eref{eq10_12} need to be modified. As before, we shall assume
that $\boldsymbol M\,{=}\,\widehat{\mathbf{A}}^{-1}\bomega +
\lambda\bgamma$, where the coefficient $\lambda$ can be found from the
condition
\[
(\boldsymbol  M+\boldsymbol  k,\bgamma)=(\bomega,\bgamma).
\]
We obtain
\[
\eqalign{
\boldsymbol  M=\widehat{\mathbf{A}}^{-1}\bomega - \big((\widehat{\mathbf{A}}^{-1}-\mathbf{E})\bomega +\boldsymbol  k,\bgamma\big)\bgamma,  \qquad
\bomega=\widehat{\mathbf{A}}(\boldsymbol M - S\bgamma),\\
S=\frac{(\widehat{\mathbf{A}}\boldsymbol  M -\boldsymbol  M -\boldsymbol  k,\bgamma)}{(\widehat{\mathbf{A}}\bgamma,\bgamma)}.
}
\]
Here $S$ coincides with the corresponding function in the bracket~\eref{eq13_} provided that $g$ is given as
\[
g=\sqrt{(\widehat{\mathbf{A}}\bgamma,\bgamma)}.
\]
In this case the Hamiltonian reads
\[
H=\frac12\bigg( (\widehat{\mathbf{A}}\boldsymbol  M,\boldsymbol  M) +
    \frac{(\widehat{\mathbf{A}}\boldsymbol  M -\boldsymbol  M -\boldsymbol  k,\bgamma)^2}{(\widehat{\mathbf{A}}\bgamma,\bgamma)} \bigg).
\]

It turns out that this modified Veselova system with gyrostatic terms still admits one additional integral of the form
$$
F_3 = (\boldsymbol  M+\boldsymbol  k,\boldsymbol  M+\boldsymbol  k).
$$

Thus, this new system is conformally Hamiltonian and integrable. Its dynamics can be further  analysed by the standard methods.

\section*{Conclusion and discussion}

We have obtained an invariant (independent of the choice of local
coordinates on $S^2$) conformally Hamiltonian representation of generalised Chaplygin
systems on $T^*S^2$ using a degenerate Poisson structure of rank 4 in the
six-dimensional space $\mathbb{R}^6(\boldsymbol M,\bgamma)$ and have shown that this
structure is a deformation of the standard Lie\,--\,Poisson bracket in
$\mathbb{R}^6(\boldsymbol M,\bgamma)$ corresponding to the Lie algebra $e(3)$.

As applications, we have considered two nonholonomic systems: the
Chaplygin ball  and Veselova problem. In this approach (after
a suitable change of parameters) they turn out to be integrable conformally
Hamiltonian systems on the same Poisson manifold with the same set of
first integrals. The above conformally Hamiltonian representation has been
generalised to the case of adding a gyrostat (although in this case there
is no analogy between these systems any more).

To the best of our knowledge, the conformally Hamiltonian description for the Veselova system with $(\bomega, \bgamma)\ne
0$ and the integrability of its gyrostatic generalisation were unknown before 
and are presented in this paper for the first time.

\medskip

This paper poses a number of questions related primarily to nonholonomic systems.

1. Can the above approach be used to obtain a conformally Hamiltonian
description for an integrable generalisation of the
Chaplygin ball rolling on a spherical base (BMF-system) found
in~\cite{bfm,bmm,tsiganov}?

2. Poisson brackets of a quite similar type are encountered in examples but with 
a Casimir function linear in $\boldsymbol M$ different from $(\boldsymbol M,\bgamma)$~\cite{biz_tsig}. It would be
interesting to find out whether such brackets can be reduced to the
standard Poisson-Lie bracket on $e^*(3)$ using the technique described above.

3. Since the Chaplygin ball problem without potential (i.e., $U(\bgamma) = 0$)
is integrable on the whole space $\mathbb{R}^6(\boldsymbol M, \bgamma)$, 
Theorem \ref{thm3} allows us to obtain a globally integrable Hamiltonian system on
$e^*(3)$, i.e., for all values of the area constant $(\boldsymbol M,\bgamma)$. As is well known, this
circumstance may be interpreted as integrability of a natural system with
a magnetic field whose additional integral is quadratic in momenta. The
issue of description of all such systems was actively discussed in the
literature. It would be interesting to interpret the system thus obtained
in the context of recent classification results by V.~Marikhin  and  V.~Sokolov
 \cite{marikhin,maryhin}.

\section*{Acknowledgments}
The work of Alexey V.\,Borisov was carried out within the framework of the
state assignment to the Udmurt State University ``Regular and Chaotic
Dynamics''. The work of Ivan S.\,Mamaev was supported by the RFBR grants
13-01-12462-ofi~m.

\end{document}